\documentclass[a4paper,12pt]{amsart}
\usepackage{amssymb}
\usepackage{ifthen}
\usepackage{graphicx}
\nonstopmode
\newtheorem{thm}{Theorem}
\newtheorem{cor}{Corollary}
\newtheorem{lem}{Lemma}

\newtheorem{claim}{Claim}

\newtheorem{conj}{Conjecture}
\newtheorem{prob}{Problem}

\theoremstyle{definition}

\newtheorem{example}{Example}

\newenvironment{rem}{%
\bigskip
\noindent \textsl{{\sl Remark. }}}{\bigskip}
\newenvironment{rems}{%
\bigskip
\noindent \textsl{{\sl Remarks. }}}{\bigskip}

\newenvironment{pf}[1][]{%
 \vskip 1mm
 \noindent
 \ifthenelse{\equal{#1}{}}%
  {{\slshape Proof. }}%
  {{\slshape #1.} }%
 }%
{\qed\medskip}

\newcounter{alphabet}
\newcounter{tmp}
\newenvironment{Thm}[1][]{\refstepcounter{alphabet}%
\bigskip%
\noindent%
{\bf Theorem \Alph{alphabet}}%
\ifthenelse{\equal{#1}{}}{}{ (#1)}%
{\bf .} \itshape}{\vskip 8pt}

\makeatletter
\newcommand{\Ref}[1]{\@ifundefined{r@#1}{}{\setcounter{tmp}{\ref{#1}}\Alph{tmp}}}
\makeatother

\newcommand{\IR}{{\mathbb R}}

\newcommand{\IC}{{\mathbb C}}
\newcommand{\ID}{{\mathbb D}}

\def\be{\begin{equation}}
\def\ee{\end{equation}}

\newcommand{\bee}{\begin{enumerate}}
\newcommand{\eee}{\end{enumerate}}

\newcommand{\blem}{\begin{lem}}
\newcommand{\elem}{\end{lem}}
\newcommand{\bthm}{\begin{thm}}
\newcommand{\ethm}{\end{thm}}
\newcommand{\bcor}{\begin{cor}}
\newcommand{\ecor}{\end{cor}}
\newcommand{\beg}{\begin{example}}
\newcommand{\eeg}{\end{example}}
\newcommand{\begs}{\begin{examples}}
\newcommand{\eegs}{\end{examples}}
\newcommand{\bdefe}{\begin{defin}}
\newcommand{\edefe}{\end{defin}}
\newcommand{\bprob}{\begin{prob}}
\newcommand{\eprob}{\end{prob}}
\newcommand{\bques}{\begin{ques}}
\newcommand{\eques}{\end{ques}}
\newcommand{\bei}{\begin{itemize}}
\newcommand{\eei}{\end{itemize}}

\newcommand{\bde}{\begin{deter}}
\newcommand{\ede}{\end{deter}}
\newcommand{\bca}{\begin{case}}
\newcommand{\eca}{\end{case}}
\newcommand{\bcl}{\begin{claim}}
\newcommand{\ecl}{\end{claim}}
\newcommand{\bcon}{\begin{conj}}
\newcommand{\econ}{\end{conj}}
\newcommand{\bcons}{\begin{conjs}}
\newcommand{\econs}{\end{conjs}}
\newcommand{\bprop}{\begin{propo}}
\newcommand{\eprop}{\end{propo}}
\newcommand{\br}{\begin{rem}}
\newcommand{\er}{\end{rem}}
\newcommand{\brs}{\begin{rems}}
\newcommand{\ers}{\end{rems}}
\newcommand{\bo}{\begin{obser}}
\newcommand{\eo}{\end{obser}}
\newcommand{\bos}{\begin{obsers}}
\newcommand{\eos}{\end{obsers}}
\newcommand{\bpf}{\begin{pf}}
\newcommand{\epf}{\end{pf}}
\newcommand{\ba}{\begin{array}}
\newcommand{\ea}{\end{array}}
\newcommand{\beq}{\begin{eqnarray}}
\newcommand{\beqq}{\begin{eqnarray*}}
\newcommand{\eeq}{\end{eqnarray}}
\newcommand{\eeqq}{\end{eqnarray*}}

\newcommand{\ds}{\displaystyle}

\newcounter{minutes}
\divide\time by 60
\newcounter{hours}
\multiply\time by 60 \addtocounter{minutes}{-\time}

\begin{document}
\title[Convolutions of slanted half-plane harmonic mappings]
{Convolutions of slanted half-plane harmonic mappings}

\thanks{
File:~\jobname .tex,
          printed: \number\day-\number\month-\number\year,
          \thehours.\ifnum\theminutes<10{0}\fi\theminutes}

\author{Liulan Li and S. Ponnusamy $^\dagger $
}
\address{Liulan Li, Department of Mathematics and Computational Science,
Hengyang Normal University, Hengyang,  Hunan 421008, People's
Republic of China}
\email{lanlimail2008@yahoo.com.cn}
\address{S. Ponnusamy, Department of Mathematics, Indian Institute of Technology Madras, Chennai 600036, India}
\email{samy@iitm.ac.in}
\keywords{Harmonic mapping, convolution, univalent, half-plane and slanted half-plane mappings, convex functions. \\
$
^\dagger$ {\tt Corresponding author}
}
\thanks{The research was supported by the Science and Technology Development Program of Hengyang (No. 2010KJ22),
NSF of Hunan (No. 10JJ4005) and Hunan Provincial Education
Department (No. Q12034)} \maketitle

\begin{abstract}
Let ${\mathcal S^0}(H_{\gamma})$ denote the class of all univalent, harmonic,  sense-preserving
and normalized mappings $f$ of the unit disk $\ID$ onto the slanted half-plane
$H_\gamma :=\{w:\,{\rm Re\,}(e^{i\gamma}w) >-1/2\}$ with an additional condition $f_{\overline{z}}(0)=0$.
Functions in this class can be constructed by the shear construction due to Clunie and Sheil-Small
which allows by examining their conformal counterpart. Unlike the conformal case, convolution
of two univalent harmonic convex mappings in $\ID$ is not necessarily even univalent in $\ID$.
In this paper, we fix $f_0\in{\mathcal S^0}(H_{0})$ and show that the
convolutions of $f_0$ and some slanted half-plane harmonic mapping are
still convex in a particular direction. The results of the paper enhance the interest
among harmonic mappings and, in particular,
solves an open problem of Dorff, et. al. \cite{DN} in a more general setting.
Finally, we present some basic examples of functions and their corresponding convolution functions
with specified dilatations, and illustrate them graphically with the help of MATHEMATICA software.
These examples explain the behaviour of the image domains.
\end{abstract}

\maketitle
\pagestyle{myheadings}
\markboth{Liulan Li and S. Ponnusamy}{Convolutions of harmonic mappings}

\section{Introduction}

In this paper, we consider the class $\mathcal H$ of complex-valued harmonic functions
$f=h+\overline{g}$ defined on the unit disk ${\mathbb D}=\{z \in {\mathbb C}:\, |z|<1\}$,
where $h$ and $g$ are analytic on $\ID$ with the form
\be\label{li1-eq1}
h(z)=z+\sum _{n=2}^{\infty}a_nz^n~\mbox{ and }~g(z)=\sum _{n=1}^{\infty}b_nz^n.
\ee
If we write $f = u + iv $, then $u$ and $v$ are real harmonic in $\ID$. Moreover,
the Jacobian of $f=h+\overline{g}$ is given by $J_f(z)=|h'(z)|^2 -|g'(z)|^2$.
Lewy's theorem implies that every harmonic function $f$ on $\ID$ is locally one-to-one
and sense-preserving on $\ID$ if and only if $J_f(z)>0$ in $\ID$. The condition
$J_f(z)> 0$ is equivalent to the existence of an analytic function $\omega $ in $\ID$
such that
\be\label{li1-eq0}
|\omega (z)|<1 ~\mbox{ for}~z\in \ID,
\ee
where $\omega (z)=g'(z)/h'(z)$ which is referred to as the \textit{complex dilatation} of $f$.
By requiring harmonic function to be sense-preserving we retain some basic properties exhibited
by analytic functions, such as the open mapping property, the argument principle, and zeros being
isolated (see \cite{DHL96}).

During the last two decades, after the publication of landmark paper of Clunie and Sheil-Small \cite{Clunie-Small-84},
the class ${\mathcal S}_H$ of sense-preserving univalent functions $f \in \mathcal H$
together with its subclasses have been extensively studied.
Let ${\mathcal S}_H^0$ be the subset of all $f\in {\mathcal S}_H$ in which $b_1=f_{\overline{z}}(0)=0$. We remark
that the familiar class ${\mathcal S}$ of normalized analytic univalent functions is contained in
${\mathcal S}_H^0$. Every  $f \in {\mathcal S}_H$ admits the complex dilatation $\omega $ of $f$
which satisfies \eqref{li1-eq0}. When $f \in {\mathcal S}_H^0$, we also have $\omega '(0)=0$.
Finally, let ${\mathcal K}_H^0$, ${\mathcal S}_H^{*0}$, and ${\mathcal C}_H^0$
denote  the subclasses of ${\mathcal S}_H^0$ mapping $\ID$ onto, respectively, convex, starlike,
and close-to-convex domains, just as
${\mathcal K}$, ${\mathcal S}^*$, and ${\mathcal C}$ are the subclasses of ${\mathcal S}$
mapping $\ID$ onto these respective domains.
The reader is referred to \cite{Clunie-Small-84,Du} for many interesting results and expositions
on planar univalent harmonic mappings.

%

Although not much is known in the literature on results about harmonic convolution of functions,
some progress has been achieved in the recent years, see \cite{Do,DN,Good02}.
We define the harmonic convolution (or Hadamard product) as follows:
For $f=h+\overline{g} \in {\mathcal H}$ with the series expansions for $h$ and $g$ as above
in \eqref{li1-eq1}, and $F = H + \overline{G}\in {\mathcal H}$, where
$$
H(z)=z+\sum_{n=2}^\infty A_nz^n~\mbox{ and }~ G(z)=\sum_{n=1}^\infty B_nz^n,
$$
we define
$$(f\ast F)(z)=z+\sum_{n=2}^\infty a_nA_nz^n+\sum_{n=1}^\infty \overline{b_nB_n}\overline{z^n}.
$$
Clearly, $f\ast F=F\ast f$.
In the case of conformal mappings, the literature about convolution theory is exhaustive.
For example, we have \cite{rs1}
$${\mathcal K}\ast {\mathcal K}\subset {\mathcal K},\quad  {\mathcal S}^*\ast {\mathcal K}\subset {\mathcal S}^*,\quad
{\mathcal C}\ast {\mathcal K}\subset {\mathcal C}
$$
settling the P\'{o}lya-Schoenberg conjecture. For some related containment relations, we refer
to \cite{samy95,PoSi96} and many other later works of Ruscheweyh.
Unfortunately, these inclusion results do not necessarily carryover to harmonic mappings.
In fact, in view of the sharp coefficient bounds for functions in ${\mathcal K}_H^0$,
if we take $f$ and $F$ in the class ${\mathcal K}_H^{*0}$, then
it will not always be true that $F\ast f\in{\mathcal K}_H^{*0}$ (does not necessary be even univalent).
On the other hand, based on the question raised by Clunie and Sheil-Small \cite{Clunie-Small-84},
several authors have studied the subclass of functions $f\in{\mathcal S}_H^0$ that map $\ID$ onto
specific domains such as horizontal strips, see \cite{HS87}. A function
$f=h+\overline{g} \in {\mathcal S}_H^0$ is called a
slanted half-plane mapping with $\gamma$ ($0\leq\gamma<2\pi$) if $f$
maps $\ID$ onto $H_\gamma :=\{w:\,{\rm Re\,}(e^{i\gamma}w) >-1/2\}$.
Using the shearing method due to Clunie and Sheil-Small \cite{Clunie-Small-84}, it
is almost easy to obtain that such a mapping has the form (see \cite[Lemma 1]{DN})
\be\label{li2-eq3}
h(z)+e^{-2i\gamma}g(z)=\frac{z}{1-e^{i\gamma}z}.
\ee
 We denote by ${\mathcal S^0}(H_{\gamma})$, the class
of all slanted half-plane mappings with $\gamma$. In the harmonic case, one can easily see that
there are infinitely many slanted half-plane mapping with a fixed $\gamma$.

For $\gamma=0$, we get the class of
right half-plane mappings $f$ that map $\ID$ onto $f(\ID)= H_0=\{w:\,{\rm Re\,}\ w >-1/2\}$
and such mappings clearly assume the form
$$h(z)+g(z)=\frac{z}{1-z}.
$$
Moreover if  $f\in {\mathcal K}_H^0$
and $\phi \in {\mathcal S^0}(H_{0})$, then $f\ast\phi$ is not
necessarily belong to ${\mathcal S}_H^0$. This can be easily seen to be true even
if  $\phi \in {\mathcal S^0}(H_{\gamma})$ for any $\gamma$.

Throughout the paper $f_0=h_0+\overline{g_0}$, where
$$ h_0(z)=\frac{z-z^2/2}{(1-z)^2} = z+\sum_{n=2}^\infty \frac{n+1}{2}z^n
=\frac{1}{2}\left (\frac{z}{1-z} +\frac{z}{(1-z)^2}\right )
$$
and
$$ g_0(z)=\frac{-z^2/2}{(1-z)^2} = -\sum_{n=2}^\infty \frac{n-1}{2}z^n
=\frac{1}{2}\left (\frac{z}{1-z} -\frac{z}{(1-z)^2}\right ).
$$
The function $f_0$, which acts as extremal for many issues concerning the convex class ${\mathcal K}_H^0$,
has the dilatation $\omega (z)=-z$. Moreover, if
$f=h+\overline{g}\in {\mathcal H}$, then the above
representation for $h_0$ and $g_0$ quickly gives that
$$f_0\ast f= h_0\ast h +\overline{g_0\ast g}= \frac{h+zh'}{2}+\overline{\frac{g-zg'}{2}}.
$$
This fact will be used in the proof of Lemma \ref{li2-th1} while determining the
dilatation of the convolution functions.
Clearly $f_0\in {\mathcal S^0}(H_{0})$, because
$h_0(z)+g_0(z)=\frac{z}{1-z}$ and  $f_0\in {\mathcal K}_H^0$ (see
\cite{Clunie-Small-84}). We observe that $f_0\ast f_0\not\in {\mathcal K}_H^0$, see
\cite[Theorem 5.7]{Clunie-Small-84}.

A domain $\Omega \subset \IC$ is said to be \textit{convex in the
direction} $\gamma$, $\gamma \in\IR$, if and only if for every $a \in \IC$, the set
$\Omega \cap \{a+te^{i\gamma}:\, t\in\IR\}$ is either connected or empty.
Dorff et. al. \cite{DN} proved

\begin{Thm}{\rm (\cite[Theorem 2]{DN})}\label{ThmA}
If $f_k\in {\mathcal S^0}(H_{\gamma_k})$, $k=1,\ 2$, and $f_1\ast
f_2$ is locally univalent in $\ID$, then $f_1\ast f_2$ is convex in
the direction $-(\gamma_1 +\gamma_2)$.
\end{Thm}

Theorem \Ref{ThmA} generalizes the result of Dorff \cite[Theorem 5]{Do} who proved it when
$\gamma _1=\gamma _2=0$. Moreover, Theorem \Ref{ThmA} implies that
for every $f\in {\mathcal S^0}(H_\gamma)$, $f\ast f_0$ is convex in
the direction $-\gamma$ provided that convolution function is locally univalent in $\ID$.
On the other hand, the following result deals with two cases for which the local
univalence of the resulting convolution function as demanded in Theorem \Ref{ThmA} is not necessary.

\begin{Thm}{\rm (\cite[Theorem 3]{DN})}\label{ThmB}
Let $f=h+\overline g\in  {\mathcal S^0}(H_0)$ with the dilatation
$\omega(z)=e^{i\theta}z^n$ $(n=1,2)$, $\theta\in\IR $.
Then $f_0\ast f\in {\mathcal S}_H^0$ and is convex in the direction of the real axis.
\end{Thm}


We now state our first result which shows that Theorem \Ref{ThmB} continues to hold in
the general setting.

\begin{thm}\label{univalent1}
Let $f=h+\overline{g}\in  {\mathcal S^0}(H_\gamma)$ with the
dilatation $\omega(z)=e^{i\theta}z^n$, where  $n=1,2$ and $\theta\in\IR$.
Then $f_0\ast f\in {\mathcal S}_H^0$ and is convex in the direction $-\gamma$.
\end{thm}

Also, we present an example when the local univalency fails for $\omega(z)=e^{i\theta}z^n$
if $n\geq 3$. We note that if $f\in {\mathcal S^0}(H_\gamma)$ with $\gamma=\pi$,
then we get the class of left half-plane mappings $f$ that map $\ID$ onto
$f(\ID)= H_\pi=\{w:\,{\rm Re\,}\ w <1/2\}$ so that
$$h(z)+g(z)=\frac{z}{1+z}.
$$
Setting $\gamma=\pi$ in Theorem \ref{univalent1} gives

\bcor
Let $f=h+\overline{g}\in {\mathcal S}^0(H_\pi)$ with the dilatation
$\omega(z)=e^{i\theta}z^n$,  where  $n=1,2$ and $\theta\in\IR$.
Then $f_0\ast f\in {\mathcal S}_H^0$.
\ecor

Recently, Bshouty and Lyzzaik \cite{BL} brought out a
collection of open problems and conjectures on planar harmonic mappings, proposed by many colleagues
throughout the past quarter of a century.  In \cite[Problem 3.26\textbf{(a)}]{BL}, Dorff et. al.
posed the following open question.

\vspace{6pt}
\noindent{\bf Problem.} \textit{
Let $f = h + \overline{g}\in {\mathcal S^0}(H_0)$ with
the dilatation $\omega(z)=(z +a)/(1+\overline{a}z)$, $(-1<a<1)$.
Then $f_0\ast f\in {\mathcal S}_H^0$ and is convex in the direction of the real
axis {\rm (\cite[Theorem 4]{DN})}. Determine other values of $a \in \mathbb{D}$ for
which the previous result holds.}

\vspace{6pt}

In \cite{LiPo1}, the present authors have solved this problem. In continuation of our investigation,
in this paper we consider this problem in a general setting
by allowing $f$ to vary in ${\mathcal S^0}(H_\gamma)$.

\begin{thm}\label{univalent2}
Let $f=h+\overline{g}\in  {\mathcal S^0}(H_\gamma)$ with the
dilatation $\omega(z)=\frac{z+a}{1+\overline{a}z}$, where
$a=|a|e^{i\theta}$, $\theta=\arg z$ and $|a|<1$. If
$$|a|^2\left(\cos^2\big(\theta-\frac{\gamma}{2}\big)+9\sin^2\big(\theta-\frac{\gamma}{2}\big)\right)\leq1
$$
but the condition
$$|a|\cos\big (\theta-\frac{\gamma}{2}\big)=-\cos\frac{3\gamma}{2}
~\mbox{ and }~ 3|a|\sin\big(\theta-\frac{\gamma}{2}\big)=-\sin\frac{3\gamma}{2}
$$
does not hold, then $f_0\ast f\in {\mathcal S}_H^0$ and is convex in the
direction $-\gamma$.
\end{thm}

Theorems \ref{univalent1} and \ref{univalent2} supplement the works of
Dorff et. al. \cite{DN} and, Clunie and Sheil-Small \cite{Clunie-Small-84}.
Finally, in Section \ref{sec3} we include important special cases of  Theorem \ref{univalent2}
which includes a solution to the Problem of Dorff et. al. \cite{BL} (see Corollary \ref{univalent4}).





\section{Main Lemmas}
\begin{lem}\label{li2-th1}
Let $f=h+\overline{g}\in  {\mathcal S^0}(H_\gamma)$ with the dilatation
$\omega (z)=g'(z)/h'(z)$. Then the dilatation $\widetilde{\omega}$ of $f_0\ast f$ is
\be\label{li2-eq8}
\widetilde{\omega}(z)=-ze^{-i\gamma}\left (\frac{\omega^2(z)+e^{2i\gamma}[\omega(z)-\frac{1}{2}z\omega'(z)]
+\frac{1}{2}e^{i\gamma}\omega'(z)}
{1+e^{-2i\gamma}[\omega(z)-\frac{1}{2}z\omega'(z)]+\frac{1}{2}e^{-i\gamma}z^2\omega'(z)}\right ).
\ee
\end{lem}
\bpf Assume the hypothesis that $f=h+\overline{g}\in  {\mathcal S^0}(H_\gamma)$ with $\omega (z)=g'(z)/h'(z)$. Then
$$g'(z)=\omega (z)h'(z)~\mbox{ and }~g''(z)=\omega'(z) h'(z)+\omega (z)h''(z).
$$
Moreover, as $f$ satisfies the condition \eqref{li2-eq3}, the first equality above gives
\be\label{li2-eq4}
h'(z)=\frac{1}{(1+e^{-2i\gamma}\omega(z))(1-e^{i\gamma}z)^2}
\ee
and therefore,
\be\label{li2-eq5} h''(z)=\frac{-(1-e^{i\gamma}z)e^{-2i\gamma}\omega'(z)
+2(1+e^{-2i\gamma}\omega(z))e^{i\gamma}}{(1+e^{-2i\gamma}\omega(z))^2(1-e^{i\gamma}z)^3}.
\ee
From the representation of $h_0$ and $g_0$, we see that
$$(h_0\ast h)(z)=\frac{h(z)+zh'(z)}{2}~\mbox{ and }~ (g_0\ast g)(z)=\frac{g(z)-zg'(z)}{2}.
$$
Therefore, as $f_0\ast f= h_0\ast h +\overline{g_0\ast g}$, we have
\be\label{li2-eq6}
\widetilde{\omega}(z)=\frac{(g_0\ast g)'(z)}{(h_0\ast h)'(z)}
=-\frac{zg''(z)}{2h'(z)+zh''(z)}= -\frac{z\omega'(z) h'(z)+\omega
(z)zh''(z)}{2h'(z)+zh''(z)}.
\ee
In view of \eqref{li2-eq4} and \eqref{li2-eq5}, after some computation \eqref{li2-eq6} takes the
desired form. \epf

The case $\gamma =0$ of Lemma \ref{li2-th1} apparently used in the proof of Theorem 3 in \cite{DN}
whereas the case $\gamma =\pi$ of Lemma \ref{li2-th1} gives

\begin{cor}\label{dilatation}
Let $f=h+\overline{g}\in  {\mathcal S^0}(H_\pi)$ with
the dilatation $\omega (z)$. Then the dilatation $\widetilde{\omega}$ of $f_0\ast f$ is given by
$$\widetilde{\omega}(z)=z\frac{\omega^2(z)+[\omega(z)-\frac{1}{2}z\omega'(z)]-\frac{1}{2}\omega'(z)}
{1+[\omega(z)-\frac{1}{2}z\omega'(z)]-\frac{1}{2}z^2\omega'(z)}.
$$
\end{cor}

The following lemma  is required for the proof of
Theorem~\ref{univalent1}.

\begin{lem}\label{dila1}
Let $f=h+\overline{g}\in  {\mathcal S^0}(H_\gamma)$ with the dilatation
$\omega(z)=g'(z)/h'(z)=e^{i\theta}z^n$ $(n=1,2$ and $\theta\in\IR)$.
Then the dilatation of $f_0\ast f$ is
\be\label{li2-eq9}
\widetilde{\omega}(z)=-z^ne^{(2\theta-\gamma)i}\left
(\frac{z^{n+1}+e^{(2\gamma-\theta)i}(1-\frac{n}{2})z
+\frac{n}{2}e^{(\gamma-\theta)i}}
{1+e^{(\theta-2\gamma)i}(1-\frac{n}{2})z^n+\frac{n}{2}e^{(\theta-\gamma)i}z^{n+1}}\right).
\ee
\end{lem}
\bpf Consider $\omega(z)=e^{i\theta}z^n$. Then $\omega'(z)=ne^{i\theta}z^{n-1}.$ Using these,
\eqref{li2-eq8} gives
\beqq
\widetilde{\omega}(z) &=&
-ze^{-i\gamma}\left (\frac{e^{2i\theta}z^{2n}+e^{2i\gamma}[e^{i\theta}z^n-\frac{1}{2}zne^{i\theta}z^{n-1}]
+\frac{1}{2}e^{i\gamma}ne^{i\theta}z^{n-1}}
{1+e^{-2i\gamma}[e^{i\theta}z^n-\frac{1}{2}zne^{i\theta}z^{n-1}]+\frac{1}{2}e^{-i\gamma}z^2ne^{i\theta}z^{n-1}}\right)
\eeqq
and a simplification  gives the desired formula \eqref{li2-eq9}.
\epf

\beg
The range of the dilatation function $\widetilde{\omega}$ in Lemma \ref{dila1} is not
contained in the unit disk $\ID$ if we assume $n\geq 3$. To see this, we choose
$\omega(z)=-z^n$. Then \eqref{li2-eq9} reduces to
$$\widetilde{\omega}(z)=-z^ne^{-i\gamma}\left
(\frac{z^{n+1}+e^{2i\gamma}(\frac{n}{2}-1)z -\frac{n}{2}e^{i\gamma}}
{1+e^{-2i\gamma}(\frac{n}{2}-1)z^n-\frac{n}{2}e^{-i\gamma}z^{n+1}}\right) =-z^ne^{-i\gamma}R(z) \quad \mbox{(say)}.
$$
It is a simple exercise to see that
$$|R(e^{i\alpha})|=1 ~\mbox{ and }~R(z)\big (\overline{R(1/{\overline{z}})}\big ) =1
$$
so that the function $R(z)$ maps the closed disk $|z|\leq 1$ onto itself and hence, $R$ can be written
as a finite Blaschke product
of order $n+1$. On the other hand, the product of the moduli of the zeros of $R$ in the unit disk $\ID$
is $n/2$. This means that there exists a point $z_0\in \ID$ such that $|\widetilde{\omega}(z_0)|>1$ if $n\geq 3$.
Thus, the restriction on $n$, namely, $n=1,2$ in Lemma \ref{dila1} becomes necessary for our investigations.
\eeg

\begin{lem}\label{dila2}
Let $f=h+\overline{g}\in  {\mathcal S^0}(H_\gamma)$ with the
dilatation $\omega(z) =\frac{z+a}{1+\overline{a}z}$, where $|a|<1$.
Then the dilatation $\widetilde{\omega}$ of $f_0\ast f$ is
$$\widetilde{\omega}(z)=
-ze^{-i(\gamma-\phi)}
\cdot \frac{(z+A)(z+B)}{(1+\overline{A}z)(1+\overline{B}z)},
$$
where $\phi =\arg ((1+\overline{a}e^{2i\gamma})/(1+ae^{-2i\gamma}))$, and
\be\label{li2-eq10}
t(z)=z^2+\frac{4a+e^{2i\gamma}(1+3|a|^2)}{2(1+\overline{a}e^{2i\gamma})}z
+\frac{2a^2+2ae^{2i\gamma}+e^{i\gamma}(1-|a|^2)}{2(1+\overline{a}e^{2i\gamma})}.
\ee
Here $-A, -B$ are the two roots of $t(z)=0$, and  $A, B$ may be equal.
$($The dilatation is well-defined provided $|A|, |B|\leq 1$ which will be discussed in the
next lemma$)$
\end{lem}\bpf
We have
$$\omega'(z)=\frac{1-|a|^2}{(1+\overline{a}z)^2} ~\mbox{ for $|a|<1$}.
$$
In view of $\omega(z) =\frac{z+a}{1+\overline{a}z}$ and the last equation, by \eqref{li2-eq8},
the dilatation $\widetilde{\omega}$ of $f_0\ast f$ takes the form
$$\widetilde{\omega}(z) = -ze^{-i\gamma}W(z),
$$
where
\beqq
W(z)& =& 
\frac{\frac{(z+a)^2}{(1+\overline{a}z)^2}+e^{2i\gamma}[\frac{z+a}{1+\overline{a}z}-\frac{1}{2}z\frac{1-|a|^2}{(1+\overline{a}z)^2}]
+\frac{1}{2}e^{i\gamma}\frac{1-|a|^2}{(1+\overline{a}z)^2}}
{1+e^{-2i\gamma}[\frac{z+a}{1+\overline{a}z}-\frac{1}{2}z\frac{1-|a|^2}{(1+\overline{a}z)^2}]
+\frac{1}{2}e^{-i\gamma}z^2\frac{1-|a|^2}{(1+\overline{a}z)^2}}
\\ &=&
\frac{2(z+a)^2+e^{2i\gamma}[2(z+a)(1+\overline{a}z)-z(1-|a|^2)]
+e^{i\gamma}(1-|a|^2)}{2(1+\overline{a}z)^2+e^{-2i\gamma}[2(z+a)(1+\overline{a}z)-z(1-|a|^2)]
+e^{-i\gamma}(1-|a|^2)z^2}
\\ &=&
\frac{2(1+\overline{a}e^{2i\gamma})z^2+[4a+e^{2i\gamma}(1+3|a|^2)]z+2a^2+2ae^{2i\gamma}+e^{i\gamma}(1-|a|^2)}
{2(1+ae^{-2i\gamma})+[4\overline{a}+e^{-2i\gamma}(1+3|a|^2)]z+[2\overline{a}^2+2\overline{a}e^{-2i\gamma}
+e^{-i\gamma}(1-|a|^2)]z^2}
\\ &=&
\left(\frac{1+\overline{a}e^{2i\gamma}}{1+ae^{-2i\gamma}}\right) \left
(\frac{z^2+\frac{4a+e^{2i\gamma}(1+3|a|^2)}{2(1+\overline{a}e^{2i\gamma})}z
+\frac{2a^2+2ae^{2i\gamma}+e^{i\gamma}(1-|a|^2)}{2(1+\overline{a}e^{2i\gamma})}}
{1+\frac{4\overline{a}+e^{-2i\gamma}(1+3|a|^2)}{2(1+ae^{-2i\gamma})}z
+\frac{2\overline{a}^2+2\overline{a}e^{-2i\gamma}+e^{-i\gamma}(1-|a|^2)}{2(1+ae^{-2i\gamma})}z^2}
\right )
\\ &=&
\left(\frac{1+\overline{a}e^{2i\gamma}}{1+ae^{-2i\gamma}}\right)\frac{t(z)}{t^*(z)}.
\eeqq
Here $t(z)$ is given by \eqref{li2-eq10} and
$$t^*(z)=z^2\big (\overline{t(1/{\overline{z}})}\big )
=1+\frac{4\overline{a}+e^{-2i\gamma}(1+3|a|^2)}{2(1+ae^{-2i\gamma})}z
+\frac{2\overline{a}^2+2\overline{a}e^{-2i\gamma}+e^{-i\gamma}(1-|a|^2)}{2(1+ae^{-2i\gamma})}z^2.
$$
Suppose that $-A, -B$ are the two roots of $t(z)=0$ ($A, B$ may be
equal). Then
$$t(z)=(z+A)(z+B)
$$
and
$$t^*(z)=z^2\,\overline{t(1/\overline{z})}=z^2\cdot\overline{(1/\overline{z}+A)(1/\overline{z}+B)}
=(1+\overline{A}z)(1+\overline{B}z).
$$
As $|(1+\overline{a}e^{2i\gamma})/(1+ae^{-2i\gamma})|=1$, the desired form for $\widetilde{\omega}(z)$
follows.
\epf

\begin{lem}\label{AB}
Let $t(z)$ be defined by \eqref{li2-eq10} so that $t(z)=(z+A)(z+B)$. Also, let
$a=|a|e^{i\theta}$, where $\theta=\arg a$ with $|a|<1$. If
\be\label{li2-eq10a}
|a|^2\left(\cos^2\big (\theta-\frac{\gamma}{2}\big )+9\sin^2\big (\theta-\frac{\gamma}{2}\big )\right)\leq1,
\ee
then $|AB|\leq1$. Moreover, $|AB|=1$ if and only if
\be\label{li2-eq11}
|a|\cos\big(\theta-\frac{\gamma}{2}\big)=-\cos\frac{3\gamma}{2}\;\ \mbox{and}\;\
3|a|\sin \big(\theta-\frac{\gamma}{2}\big)=-\sin\frac{3\gamma}{2}.
\ee
\end{lem}
\bpf By the definition of $t(z)=(z+A)(z+B)$, it is clear that
$$AB=\frac{2a^2+2ae^{2i\gamma}+e^{i\gamma}(1-|a|^2)}{2(1+\overline{a}e^{2i\gamma})}
=\frac{2a(a+e^{2i\gamma})+e^{i\gamma}(1-|a|^2)}{2(1+\overline{a}e^{2i\gamma})}.
$$
We look for a condition on $a\in \ID$ such that $|AB|\leq1$. Now, a computation leads to
$$\left|2a(a+e^{2i\gamma})+e^{i\gamma}(1-|a|^2)\right|^2-4\left|1+\overline{a}e^{2i\gamma}\right|^2=(1-|a|^2)v(a),
$$
where $v(a)$ is real and
$$
v(a)=4{\rm Re\,}\big (a^2e^{-i\gamma}\big )+ 4{\rm Re\,}\big (a(e^{i\gamma}-2e^{-2i\gamma})\big )-3-5|a|^2.
$$
Now, we let $a=|a|e^{i\theta}$. Then, $v(a)$ reduces to
\be\label{li2-eq13}
v(a)=4|a|^2\cos(2\theta-\gamma)+4|a|\cos(\theta+\gamma)-8|a|\cos(\theta-2\gamma)-3-5|a|^2.
\ee
If we use the cosine doubling formula $\cos 2\phi =1-2\sin^2\phi$, and then replace $\theta+\gamma$ and $\theta-2\gamma$
respectively by $\theta -\frac{\gamma}{2}+\frac{3\gamma}{2}$ and $\theta-\gamma -\frac{3\gamma}{2}$,
by a simplification, $v(a)$ takes the form
$$
v(a)=-\left(|a|\cos\big (\theta-\frac{\gamma}{2}\big )+2\cos\frac{3\gamma}{2}\right)^2
-\left(3|a|\sin\big (\theta-\frac{\gamma}{2}\big )+2\sin\frac{3\gamma}{2}\right)^2+1.
$$

By  \eqref{li2-eq10a}, we observe that
$P_1\left(|a|\cos(\theta-\frac{\gamma}{2}),3|a|\sin(\theta-\frac{\gamma}{2})\right)$
is a point that lies on the closed disk $|z|\leq 1$ whereas the point
$P_2\left(-2\cos\frac{3\gamma}{2},-2\sin\frac{3\gamma}{2}\right)$ lies on the circle $|z|=2$.
Thus, the distance between the points $P_1$ and $P_2$ must be at least $1$. That is,
$$\sqrt{\left (|a|\cos\big (\theta-\frac{\gamma}{2}\big)+2\cos\frac{3\gamma}{2}\right )^2+
\left (3|a|\sin \big(\theta-\frac{\gamma}{2}\big) +2\sin\frac{3\gamma}{2}\right )^2}\geq1
$$
which is equivalent to saying that $v(a)\leq 0$, i.e $|AB|\leq 1$. Moreover, in the above inequality,
equality holds if and only if the point $P_1$ is the middle point of the line segment joining
$P_2$ and the origin. This gives the condition \eqref{li2-eq11}. In other words, if \eqref{li2-eq10a} holds
but not the \eqref{li2-eq11}, then $v(a)<0$ and hence, strict inequality $|AB|<1$ holds. If \eqref{li2-eq11}
holds, then $v(a)=0$ and hence, $|AB|=1$.  The proof is complete.
\epf

\section{Proofs of Main Theorems and their consequences}\label{sec3}

\subsection*{Proof of Theorem \ref{univalent1}}
In view of Theorem \Ref{ThmA}, it suffices to show that $f_1\ast f_2$ is locally univalent in $\ID$.
To prove this, first we consider the case $n=1$ so that $\omega(z)=e^{i\theta}z$.
Then, by the formula \eqref{li2-eq9}, the dilatation $\widetilde{\omega}$ of $f_0\ast f$ becomes
$$\widetilde{\omega}(z)=-ze^{(2\theta-\gamma)i}\left
(\frac{z^{2}+\frac{1}{2}e^{(2\gamma-\theta)i}z
+\frac{1}{2}e^{(\gamma-\theta)i}}
{1+\frac{1}{2}e^{(\theta-2\gamma)i}z+\frac{1}{2}e^{(\theta-\gamma)i}z^{2}}\right)
= -ze^{(2\theta-\gamma)i}\frac{t(z)}{t^*(z)},
$$
where $t(z)=z^{2}+\frac{1}{2}e^{(2\gamma-\theta)i}z +\frac{1}{2}e^{(\gamma-\theta)i}$ and
$$t^*(z)=z^2\, \overline{t{(1/\overline{z})}}=1+\frac{1}{2}e^{(\theta-2\gamma)i}z+\frac{1}{2}e^{(\theta-\gamma)i}z^{2}.
$$
Clearly if $z_0$ is a zero of $t(z)$, then $1/\overline{z_0}$ is a zero of $t^*(z)$.
Therefore, we may  write the last expression as
$$\widetilde{\omega}(z)=-ze^{(2\theta-\gamma)i}\frac{(z+A)(z+B)}{(1+\overline{A} z)(1+\overline{B} z)}.
$$
We observe that $A$ and $B$ are nonzero complex numbers such that
$$A+B =\frac{1}{2}e^{(2\gamma-\theta)i}~ \mbox{ and }~ AB =\frac{1}{2}e^{(\gamma-\theta)i}.
$$
It is easy to see that $A, B\in \overline{\ID}$. Observe that $|\frac{1}{2}e^{(\gamma-\theta)i}|<1$ and
the only zero of
$$\frac{t(z)-\frac{1}{2}e^{(\gamma-\theta)i}t^*(z)}{z}=\frac{3}{4}z
+\frac{1}{2}e^{(2\gamma-\theta)i}-\frac{1}{4}e^{-i\gamma},
$$
namely, $\frac{1}{3}e^{-i\gamma}-\frac{2}{3}e^{(2\gamma-\theta)i}$, clearly lies
in $\overline{\ID}$. According to Cohn's Rule (\cite{Co22} or see \cite{RhSch02}),
the two zeros of $t(z)$, namely $-A$ and $-B$, must lie in $\overline{\ID}$.
This observation gives that $|\widetilde{\omega}(z)|<1$ in $\ID$.

Next, we consider the case $n=2$ so that $\omega(z)=e^{i\theta}z^2$. In this case,
the formula \eqref{li2-eq9} takes the form
$$\widetilde{\omega}(z)=-z^2e^{(2\theta-\gamma)i}
\left(\frac{z^3+e^{(\gamma-\theta)i}}{1+e^{-(\gamma-\theta)i}z^3}\right),
$$
which clearly implies that $|\widetilde{\omega}(z)|<1$ for $z\in \ID$.
According to Lewy's  theorem, it turns outs that $f_1\ast f_2$ is locally univalent in $\ID$
and hence, the desired conclusion follows from Theorem \Ref{ThmA}.
\qed

\subsection*{Proof of Theorem \ref{univalent2}}
By Lemma \ref{AB} and the hypothesis, we have
%
$|AB|<1$. Then at least one of $A, B$ is in $\ID$.
Without loss of generality, we may assume that $A\in\ID$.
Next, consider the function $t(z)$ defined by \eqref{li2-eq10} in the form
$$t(z)=z^2+a_1z+a_0=(z+A)(z+B).
$$
Then the function
$$t_1(z)=\frac{t(z)-a_0t^*(z)}{z}=(1-|a_0|^2)z+a_1-a_0\overline{a_1}
$$
has a zero at
\be\label{li2-eq5a}
z_0= \frac{a_0\overline{a_1}-a_1}{1-|a_0|^2} =\frac{A(|B|^2-1)+B(|A|^2-1)}{1-|AB|^2},
\ee
which, after simplification, is equivalent to
\be\label{li2-eq6a}
z_0=e^{-i\gamma}\frac{6a^2e^{-i\gamma}+8ae^{i\gamma}-4\overline{a}e^{2i\gamma}-3|a|^2+2e^{3i\gamma}-1}
{4{\rm Re\,}\big (a^2e^{-i\gamma}\big )+ 4{\rm Re\,}\big (a(e^{i\gamma}-2e^{-2i\gamma})\big )-3-5|a|^2}.
\ee

\noindent{\bf Claim \textbf{(a).}}\quad \textit{$B\in\overline{\ID}$
if and only if $|z_0|\leq1$.}

\vspace{4pt}

By a routine computation, it can be easily seen that
$$|A(|B|^2-1)+B(|A|^2-1)|^2-(1-|AB|^2)^2=-(1-|A|^2)(1-|B|^2)|1-A\overline{B}|^2.$$
As $|AB|< 1$ and $|A|<1$, the last equation and \eqref{li2-eq5a}
show that Claim \textbf{(a)} holds. Indeed $|B|<1$ if and only if $|z_0|<1$, and $|B|=1$ if and only if
$|z_0|=1$.

\vspace{4pt}

\noindent{\bf Claim \textbf{(b).}}\quad \textit{$|z_0|\leq1$ if and
only if \eqref{li2-eq10a} holds.}

\vspace{4pt}
We may conveniently write $z_0$ as
$$z_0=e^{-i\gamma}\frac{u(a)}{v(a)}
$$
where $v(a)$ is defined by \eqref{li2-eq13} and
$$u(a)=6a^2e^{-i\gamma}+8ae^{i\gamma}-4\overline{a}e^{2i\gamma}-3|a|^2+2e^{3i\gamma}-1.
$$
We observe that $|z_0|\leq 1$ if and only if $|u(a)|\leq |v(a)|$, i.e. $|u(a)|^2-|v(a)|^2\leq 0$. In order
to deal with the later inequality, we consider
\be\label{li2-eq15}
|u(a)|^2-|v(a)|^2= ({\rm Im\,}u(a))^2 +({\rm Re\,}u(a)-v(a))({\rm Re\,}u(a) +v(a))
\ee
and each term needs to be simplified. First we find that
$$ {\rm Re\,}u(a)=6|a|^2\cos(2\theta-\gamma)+8|a|\cos(\theta+\gamma)-4|a|\cos(\theta-2\gamma)-3|a|^2+2\cos3\gamma-1,
$$
and
\beqq
{\rm Im\,}u(a) &= &6|a|^2\sin(2\theta-\gamma)+8|a|\sin(\theta+\gamma)+4|a|\sin(\theta-2\gamma)+2\sin3\gamma \\
&=& 12|a|^2\sin\big (\theta-\frac{\gamma}{2}\big)\cos\big(\theta-\frac{\gamma}{2}\big)
+12|a|\sin \big(\theta-\frac{\gamma}{2}\big)\cos\frac{3\gamma}{2} \\
&& \hspace{1cm} +4|a|\cos \big(\theta-\frac{\gamma}{2}\big)\sin\frac{3\gamma}{2}
+4\sin\frac{3\gamma}{2}\cos\frac{3\gamma}{2}\\
&=&4 \left(|a|\cos \big(\theta-\frac{\gamma}{2}\big)+\cos\frac{3\gamma}{2}\right)
\left( 3|a|\sin\big(\theta-\frac{\gamma}{2}\big)+\sin\frac{3\gamma}{2}\right).
\eeqq
Also, we see that
\beqq
{\rm Re\,}u(a)-v(a) &=& 2|a|^2(\cos(2\theta-\gamma)+1) +4|a|[\cos(\theta+\gamma)+\cos(\theta-2\gamma)]
+2(\cos3\gamma +1)\\
&=& 4|a|^2\cos^2\big(\theta-\frac{\gamma}{2}\big) +8|a|
\cos\big(\theta-\frac{\gamma}{2}\big)\cos\frac{3\gamma}{2}+4\cos^2\frac{3\gamma}{2}\\
&=& 4\left(|a|\cos \big(\theta-\frac{\gamma}{2}\big)+\cos\frac{3\gamma}{2}\right) ^2
\eeqq
and similarly,
\beqq
{\rm Re\,}u(a)+v(a) &=& 2|a|^2(5\cos(2\theta-\gamma)-4) +12|a|[\cos(\theta+\gamma)-\cos(\theta-2\gamma)]
+2(\cos3\gamma -2)\\
&=& 2\left [|a|^2\Big(1-10\sin^2\big(\theta-\frac{\gamma}{2}\big)\Big ) -12|a|
\sin\big(\theta-\frac{\gamma}{2}\big)\sin\frac{3\gamma}{2} -2\sin^2\frac{3\gamma}{2}-1
\right ].
\eeqq
Using the above expressions, \eqref{li2-eq15} takes the form
\beqq
|u(a)|^2-|v(a)|^2& =& 8\left[|a|\cos\big(\theta-\frac{\gamma}{2}\big)
+\cos\frac{3\gamma}{2}\right]^2\left[2 \left( 3|a|\sin\big(\theta-\frac{\gamma}{2}\big)+\sin\frac{3\gamma}{2}
\right)^2 \right .\\
&& \left .\hspace{.2cm} +  |a|^2\Big(1-10\sin^2\big(\theta-\frac{\gamma}{2}\big)\Big ) -12|a|
\sin\big(\theta-\frac{\gamma}{2}\big)\sin\frac{3\gamma}{2} -2\sin^2\frac{3\gamma}{2}-1\right]\\
&=&  8\left[|a|\cos\big(\theta-\frac{\gamma}{2}\big)
+\cos\frac{3\gamma}{2}\right]^2\left[\left(\cos^2 \big(\theta-\frac{\gamma}{2}\big)
+9\sin^2\big(\theta-\frac{\gamma}{2}\big)\right)|a|^2-1\right].
\eeqq
Since $|a|<1$, the last equality shows that $|z_0|<1$ if and
only if
$$|a|^2\left(\cos^2(\theta-\frac{\gamma}{2})+9\sin^2(\theta-\frac{\gamma}{2})\right)<1.
$$
Also, $|z_0|=1$ if and only if
$$|a|^2\left(\cos^2\big(\theta-\frac{\gamma}{2}\big)+9\sin^2\big(\theta-\frac{\gamma}{2}\big)\right)
=1.
$$

In conclusion, the assumption and Claims \textbf{(a)} and
\textbf{(b)} imply that $A\in\ID$ and $B\in\overline{\ID}$. We
obtain that $|\widetilde{\omega}(z)|<1$ for each $z\in\ID$. Thus, by
Theorem \Ref{ThmA}, we deduce that $f_0\ast f\in {\mathcal S}_H^0$
and $f_0\ast f$ is convex in the direction $-\gamma$.
\qed



\begin{cor}\label{univalent3}
Let $f=h+\overline{g}\in  {\mathcal S^0}(H_\gamma)$ with the
dilatation $\omega(z)=\frac{z+a}{1+\overline{a}z}$, where
$a=|a|e^{i\theta}$, $\theta=\arg z$, $|a|<1$ and $\gamma\in\{0,
\frac{2}{3}\pi, \frac{4}{3}\pi\}$. If \eqref{li2-eq10a} holds,
then $f_0\ast f\in {\mathcal S}_H^0$ and is convex in the direction
$-\gamma$.
\end{cor}
\bpf
If $\gamma\in\big \{0, \frac{2}{3}\pi, \frac{4}{3}\pi\big \}$ in Lemma \ref{AB},  then obviously
$$\Big|\cos\frac{3\gamma}{2}\Big|=1
$$
and so the first equality in \eqref{li2-eq11} holds if and only if $|a|=1$, which contradicts
the fact that $|a|<1$.  In other words, if $\gamma\in\big \{0, \frac{2}{3}\pi, \frac{4}{3}\pi\big \}$
then $|AB|\neq 1$ and so in this case, the strict inequality $|AB|<1$ holds under
the condition \eqref{li2-eq10a}. This completes the proof.
\epf

The case $\gamma =0$ of Corollary \ref{univalent3} gives the following result (see \cite[Theorem]{LiPo1}).

\begin{cor}\label{univalent4}
Let $f=h+\overline{g}\in  {\mathcal S^0}(H_0)$ with the dilatation $\omega(z)=\frac{z+a}{1+\overline{a}z}$
$(|a|<1)$. If
$$({\rm Re\,}a)^2+9({\rm Im\,}a)^2\leq 1,
$$
then $f_0\ast f\in {\mathcal S}_H^0$ and is convex in the direction of the real axis.
\end{cor}

Finally, the case $\gamma =\pi $ of Theorem \ref{univalent2} gives

\bcor
Let $f=h+\overline{g}\in {\mathcal S}_H^0$ with $h(z)+g(z)=\frac{z}{1+z}$ and the dilatation
$\omega(z)=\frac{z+a}{1+\overline{a}z}$, where  $|a|<1$ with ${\rm Im\,}a\neq 0$. If
$$9({\rm Re\,}a)^2+({\rm Im\,}a)^2\leq 1,
$$
then $f_0\ast f\in {\mathcal S}_H^0$ and is convex in the direction of the real axis.
\ecor

\section{Some examples}

\begin{example}\label{eg1}
Let $f=h+\overline{g}\in {\mathcal S^0}(H_\gamma)$ with
$\gamma=\pi/2$ and the dilatation $\omega(z)=z$. Then
$$h(z)-g(z)=\frac{z}{1-iz} ~\mbox{ and }~g'(z)=zh'(z).
$$
Solving these yield
$$h'(z)=\frac{1}{(1-z)(1-iz)^2}=\frac{1-i}{2}\cdot \frac{1}{(1-iz)^2}+ \frac{1}{2}\left (\frac{1}{1-iz}+\frac{i}{1-z}
\right ).$$
Integration from $0$ to $z$ gives
$$h(z)=\frac{i}{2}\log \left ( \frac{1-iz}{1-z}\right )+\frac{1-i}{2}\frac{z}{1-iz}
$$
and therefore,
$$g(z)
=\frac{i}{2}\log\left ( \frac{1-iz}{1-z}\right )-\frac{1+i}{2}\frac{z}{1-iz}.
$$
By the convolution, we have
$$f_0\ast f=h_0\ast h+ \overline{g_0\ast g}
$$
so that
$$h_0\ast h=\frac{h+zh'}{2}=\frac{i}{4}\log
\frac{1-iz}{1-z}+\frac{1-i}{4}\frac{z}{1-iz}+\frac{z}{2(1-z)(1-iz)^2}
$$
and
$$g_0\ast g=\frac{g-zg'}{2}=\frac{i}{4}\log
\frac{1-iz}{1-z}-\frac{1+i}{4}\frac{z}{1-iz}-\frac{z^2}{2(1-z)(1-iz)^2}.
$$
Applying Lemma \ref{li2-th1} with $\gamma=\pi/2$, we get
$$\widetilde{\omega}(z)=iz\left (\frac{z^{2}-\frac{1}{2}z +\frac{i}{2}}
{1-\frac{1}{2}z-\frac{i}{2}z^{2}}\right).
$$
The real and imaginary parts of $f_0\ast f=h_0\ast h+ \overline{g_0\ast g}$
may be written explicitly as follows:
$${\rm Re\,}(f_0\ast f)={\rm Re\,}(h_0\ast h+ g_0\ast g)={\rm
Re\,}\left(\frac{i}{2}\log
\left ( \frac{1-iz}{1-z}\right )+\frac{z(1-i-z)}{2(1-iz)^2} \right),
$$
and
$${\rm Im\,}(f_0\ast f)={\rm Im\,}(h_0\ast h- g_0\ast g)={\rm Im\,}\left(\frac{z}{2(1-iz)}
+\frac{z+z^2}{2(1-z)(1-iz)^2} \right).
$$
The images of the unit disk $\ID$ under $f$ and $f_0\ast f$ are shown in Figure \ref{fig1}
as plots of the images of equally spaced radial segments and concentric circles.
Similar comments apply for the other two figures in the next two examples.
\begin{figure}[htp]
\begin{center}
\includegraphics[width=7.5cm]{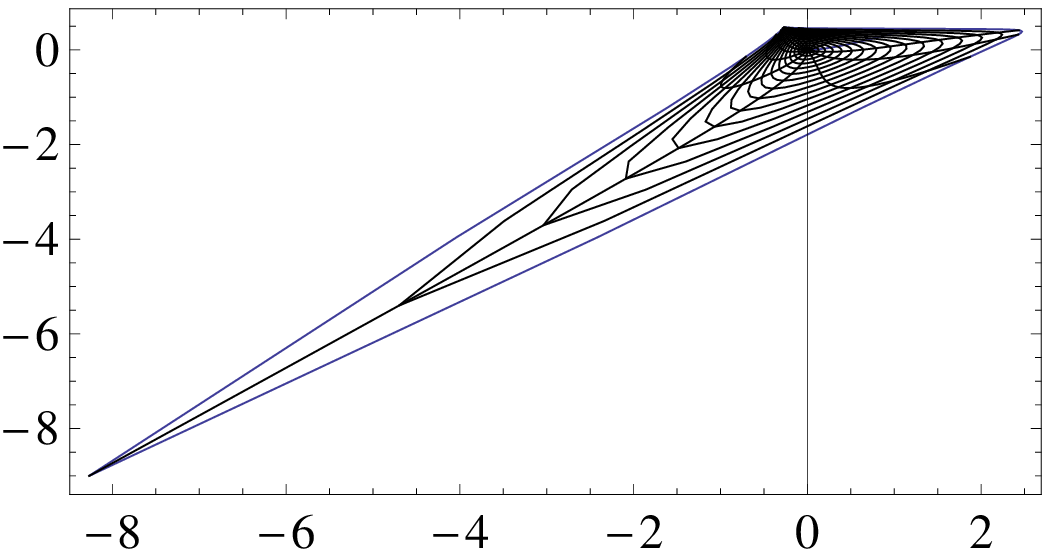}
\hspace{0.5cm}
\includegraphics[width=7.5cm]{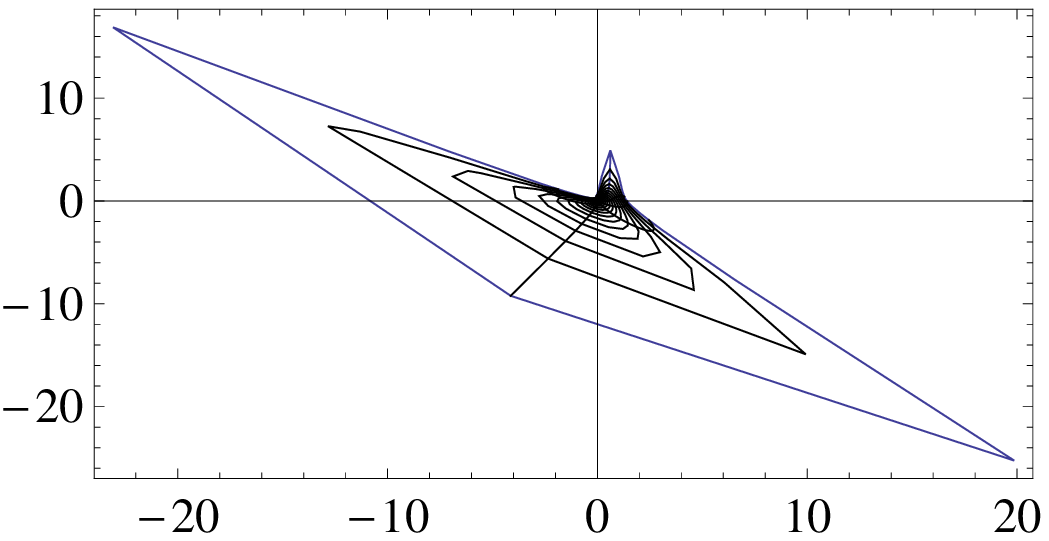}
\end{center}
\hspace{1.5cm}$f(z)$ \hspace{4.5cm}  $F(z)$
\caption{Images of $f$ and $f_0\ast f$\label{fig1}}
\end{figure}

If we need to exactly describe the image domains, one may proceed by introducing
$$\zeta =\frac{1-iz}{1-z}=\xi+i\eta  ~~(\xi>\eta ) ,
$$
so that $z=(\zeta -1)/(\zeta -i). $
Applying this transformation with $\zeta =\xi+i\eta,$ we get
$${\rm Re\,}(f_0\ast f)={\rm  Re}\left(\frac{i}{2}\log
\zeta +\frac{\zeta ^2-1}{4\zeta ^2}\right)
=-\frac{1}{2}\arctan\frac{\eta}{\xi}+\frac{1}{4}-\frac{\xi^2-\eta^2}{4(\xi^2+\eta^2)^2},$$
and
\beqq
{\rm Im\,}(f_0\ast f)& =&
{\rm Im }\left(\frac{\zeta -1}{2(1-i)\zeta }+\frac{(\zeta -1)(\zeta -i)(2\zeta -1-i)}{2(1-i)^3\zeta ^2}
\right)\\
&=&
\frac{1}{4}+\frac{\xi-\eta}{4}-\frac{3(\xi-\eta)}{4(\xi^2+\eta^2)}+\frac{\xi^2-\eta^2}{4(\xi^2+\eta^2)^2}.
\eeqq
A careful analysis may be done in order to explain the image domain of the convolution function. Here we avoid
the computation although we just would like indicate the procedure for discussion.
\end{example}

\begin{example}\label{eg2}
Let $f=h+\overline{g}\in {\mathcal S^0}(H_\gamma)$ with $\gamma=\pi$ and the dilatation $\omega(z)=z$.
Then
$$h(z)+g(z)=\frac{z}{1+z},\;\ h'(z)=\frac{1}{(1+z)^3},\;\ g'(z)=\frac{z}{(1+z)^3},$$
and therefore,
$$h(z)=\frac{z^2+2z}{2(1+z)^2},
$$
and
$$g(z)
=\frac{z}{1+z}-\frac{z^2+2z}{2(1+z)^2}=\frac{z^2}{2(1+z)^2}.
$$
As before, we easily have $f_0\ast f=h_0\ast h+ \overline{g_0\ast g}$ with
$$h_0\ast h=\frac{h+zh'}{2}
=\frac{z(z^2+3z+4)}{4(1+z)^3}
$$
and
$$g_0\ast g=\frac{g-zg'}{2}
=-\frac{z^2(1-z)}{4(1+z)^3}.
$$
The dilatation $\widetilde{\omega}$ of $f_0\ast f$ is
$$\widetilde{\omega}(z)=z\left (\frac{z^{2}+\frac{1}{2}z -\frac{1}{2}}
{1+\frac{1}{2}z-\frac{1}{2}z^{2}}\right).
$$
Further, the real and imaginary parts of $f_0\ast f=h_0\ast h+
\overline{g_0\ast g}$ are given by
$${\rm Re\,}(f_0\ast f)={\rm Re\,}(h_0\ast h+ g_0\ast g)={\rm
Re}\left(\frac{z(z^2+z+2)}{2(1+z)^3} \right)
$$
and
$${\rm Im\,}(f_0\ast f)={\rm Im\,}(h_0\ast h- g_0\ast g)={\rm Im\,}\frac{z}{(1+z)^2},
$$
respectively.
The images of the unit disk $\ID$ under $f$ and $f_0\ast f$ are shown in Figure \ref{fig2}.

\begin{figure}[htp]
\begin{center}
\includegraphics[width=7.5cm]{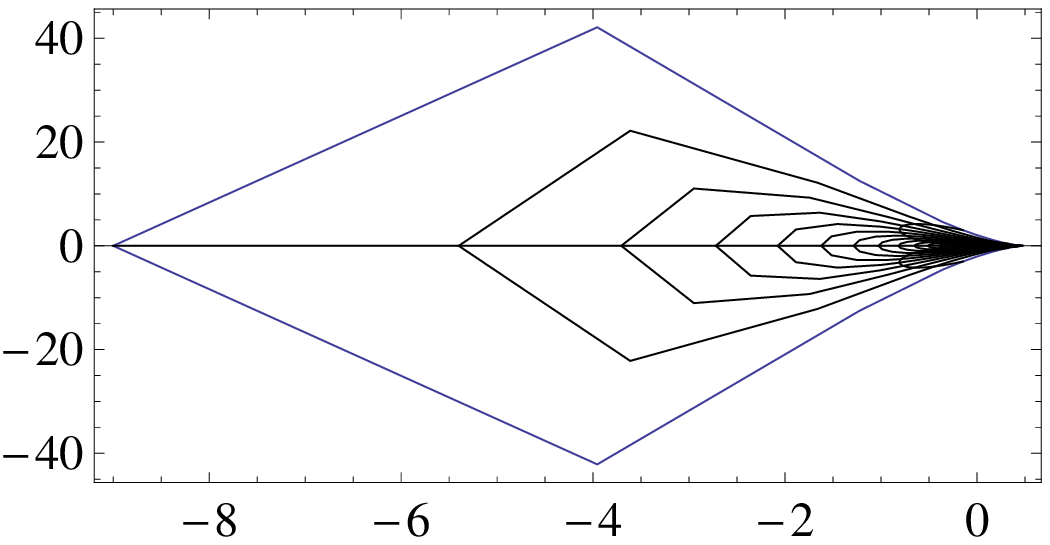}
\hspace{0.5cm}
\includegraphics[width=7.5cm]{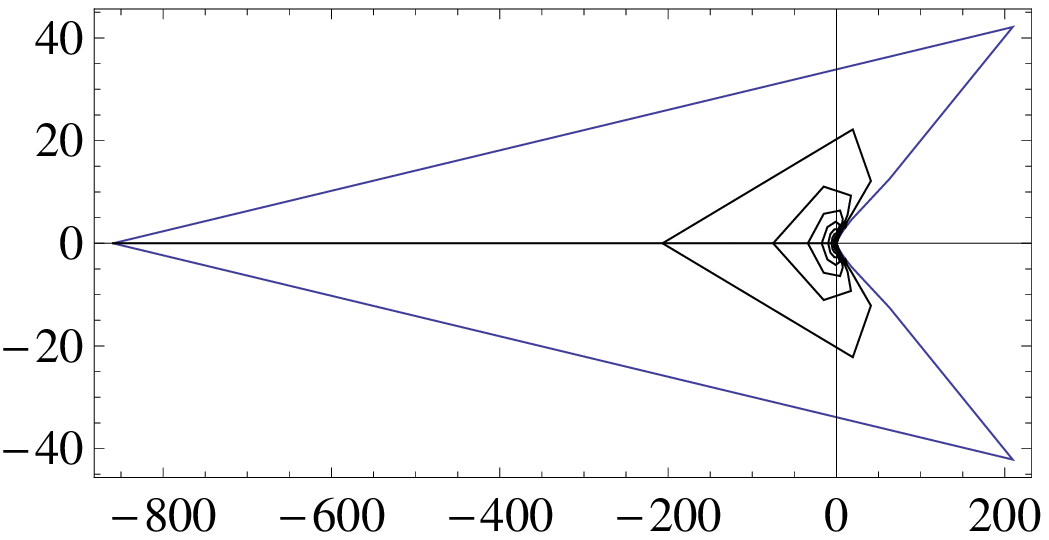}
\end{center}
\hspace{1.5cm}$f(z)$ \hspace{4.5cm}  $F(z)$
\caption{Images of $f$ and $f_0\ast f$\label{fig2}}
\end{figure}

If we let
$$\zeta =\frac{1-z}{1+z}=\xi+i\eta ~~(\xi>0 ),
$$
then $z=(1-\zeta)/(1+\zeta )$ so that
$${\rm Re\,}(f_0\ast f)=\frac{1}{8}{\rm Re}\left(-\zeta ^3-\zeta+2
  \right)=-\frac{\xi^3}{8}+\frac{3}{8}\xi\eta^2-\frac{\xi}{8}+\frac{1}{4}
$$
and
$${\rm Im\,}(f_0\ast f)={\rm Im}\left(\frac{1-\zeta ^2}{4}\right)=-\frac{\xi\eta}{2}.
$$
Again, these observations suffice for the discussion of the image domain of $f_0\ast f$.
\end{example}

\begin{example}\label{eg3}
Let $f=h+\overline{g}$ be a harmonic mapping of $\ID$ such that
$$h(z)+g(z)=\frac{z}{1+z} ~\mbox{ and }~\omega(z)=-z^2.
$$
Then $g'(z)=-z^2h'(z)$ and as before, we see that
$$ h'(z)=\frac{1}{(1-z)(1+z)^3} ~\mbox{ and }~ g'(z)=-\frac{z^2}{(1-z)(1+z)^3}.
$$
Integration gives
$$h(z)=\frac{1}{8}\log \left ( \frac{1+z}{1-z}\right ) + \frac{1}{4}\left (\frac{z}{1+z}- \frac{1}{(1+z)^2}+1\right ),
$$
$$g(z)=-\frac{1}{8}\log \left ( \frac{1+z}{1-z}\right ) + \frac{1}{4}\left (\frac{3z}{1+z}+ \frac{1}{(1+z)^2}-1\right ),
$$
and $f(\ID)= \{w:\,{\rm Re\,}\ w <1/2\}$. We observe that
$$f(e^{i\theta})
=\left \{
\begin{array}{ll}
\ds \frac{1}{2} +\left (\frac{\pi}{16} +\frac{1}{4}\tan \left (\frac{\theta}{2}\right) \right )
& \mbox{ if $0<\theta <\pi$}\\
\ds \frac{1}{2} +\left (-\frac{\pi}{16} +\frac{1}{4}\tan \left (\frac{\theta}{2}\right)  \right )
& \mbox{ if $\pi <\theta <2\pi$}.
\end{array}
\right .
$$
Next we note that $f_0\ast f=h_0\ast h+ \overline{g_0\ast g}$ with
$$h_0\ast h=\frac{1}{16}\log \left ( \frac{1+z}{1-z}\right )
 + \frac{1}{8}\left (\frac{z}{1+z}-\frac{1}{(1+z)^2}+1\right )+
\frac{z}{2(1-z)(1+z)^3}
$$
and
$$g_0\ast g=-\frac{1}{16}\log \left ( \frac{1+z}{1-z}\right )
 + \frac{1}{8}\left (\frac{3z}{1+z}+\frac{1}{(1+z)^2}-1\right )+
\frac{z^3}{2(1-z)(1+z)^3}.
$$
It is also easy to see that  $\widetilde{\omega}$ of $f_0\ast f$ is $\widetilde{\omega}(z)=z^2$.
The images of the unit disk $\ID$ under $f$ and $f_0\ast f$ are shown in Figure \ref{fig3}.

\begin{figure}[htp]
\begin{center}
\includegraphics[width=7.5cm]{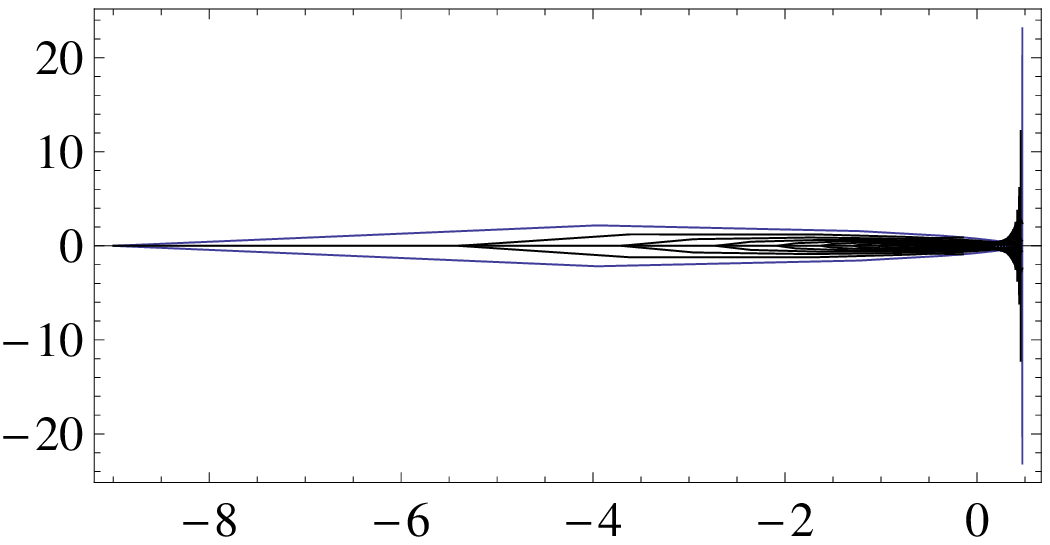}
\hspace{0.5cm}
\includegraphics[width=7.5cm]{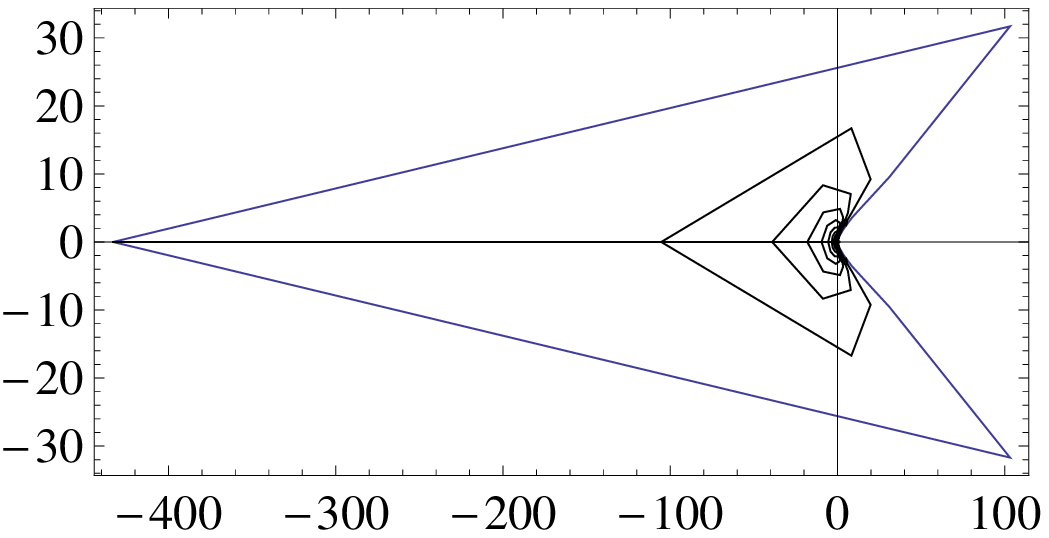}
\end{center}
\hspace{1.5cm}$f(z)$ \hspace{4.5cm}  $F(z)$
\caption{Images of $f$ and $f_0\ast f$\label{fig3}}
\end{figure}

Finally, if we let
$$\zeta =\frac{1+z}{1-z}=\xi+i\eta ~~(\xi>0),
$$
then $z= (\zeta-1)/(\zeta+1 )$ so that
\beqq
{\rm Re\,}(f_0\ast f) & =& \frac{1}{2}{\rm
Re}\left(\frac{z}{1+z}+\frac{z(1+z^2)}{2(1-z)(1+z)^3}\right) \\
& =& \frac{1}{16}{\rm Re}\left(\zeta+4 -\frac{4}{\zeta}-\frac{1}{\zeta^3} \right)\\
&=& \frac{\xi}{16}+\frac{1}{4}-\frac{\xi}{4(\xi^2+\eta^2)}-\frac{\xi^3-3\xi\eta^2}{16(\xi^2+\eta^2)^3},
\eeqq
and
\beqq
{\rm Im\,}(f_0\ast f) &= &{\rm Im}\left(\frac{1}{8}\log \left ( \frac{1+z}{1-z}\right )
-\frac{z}{4(1+z)}-\frac{1}{4(1+z)^2}+\frac{z}{2(1+z)^2} \right)\\
&=&{\rm Im}\left(\frac{1}{8}\log\zeta-\frac{3}{16\zeta^2}\right)\\
&=&\frac{1}{8}\arctan\frac{\eta}{\xi}+\frac{3\xi\eta}{8(\xi^2+\eta^2)^2}.
\eeqq
These observations help to analyze the image of $\ID$ under $f_0\ast f$ and we avoid the details.
\end{example}


\begin{thebibliography}{99}

\bibitem{BL} D.~Bshouty and A.~Lyzzaik, Problems and conjectures in planar harmonic mappings,
In: Proc. ICM2010 Satellite Conf. Int. Workshop on Harmonic and Quasiconformal
Mappings, IIT Madras, Aug. 09-17, 2010, ed. by D. Minda, S. Ponnusamy,
and N. Shanmugalingam,
\textit{J. Analysis} {\bf 18}(2010), 69--81.
See also, {\tt http://mat.iitm.ac.in/home/samy/publi${\rm c_{-}}$html/Journal.htm}

\bibitem{Clunie-Small-84}
J.~G.~Clunie and T.~Sheil-Small, Harmonic univalent functions,
\textit{Ann. Acad. Sci. Fenn. Ser. A.I.} {\bf 9} (1984), 3--25.

\bibitem{Co22} A.~Cohn,
\"{U}ber die Anzahl der Wurzeln einer algebraischen  Gleichung in einem Kreise,
\textit{Math. Z.} {\bf 14}(1922), 110--148.

\bibitem{Do} M.~Dorff,
Convolutions of planar harmonic convex mappings,
\textit{Comp. Vari. Theo. Appl.} {\bf 45}(2001), 263--271.

\bibitem{DN} M.~Dorff, M.~Nowak and M.~Wo{\l}oszkiewicz,
Convolutions of harmonic convex mappings,
\textit{Comp. Vari. Elliptic Eqn.} (2011),
DOI: 10.1080/17476933.2010.487211; arXiv:0903.1595.

\bibitem{Du} P.~Duren,
Harmonic mappings in the plane,
Cambridge Tracts in Mathematics, 156. Cambridge Univ. Press, Cambridge, 2004.

\bibitem{DHL96} P.~Duren,  W.~Hengartner and R. Laugesen,
The augment principle for harmonic mappings,
\textit{Amer. Math. Monthly} \textbf{103}(5)(1996), 411--415.

\bibitem{Good02} M.~R.~Goodloe,
Hadamard products of convex harmonic mappings,
\textit{Comp. Vari. Theo. Appl.} {\bf 47}(2)(2002), 81--92.

\bibitem{HS87} W.~Hengartner and G.~Schober,
Univalent harmonic functions,
\textit{Trans. Amer. Math. Soc.} {\bf 299}(1)(1987), 1-–31.

\bibitem{LiPo1} Liulan Li and S.~Ponnusamy, Solution to an open problem on convolutions
of harmonic mappings,
\textit{Comp. Vari. Elliptic Eqn.} (2012), Accepted.

\bibitem{samy95} S.~Ponnusamy,
P\'{o}lya-Schoenberg conjecture for Carath\'{e}odory functions,
\textit{J. London Math. Soc.} {\bf 51}(2)(1995), 93--104.

\bibitem{PoSi96} S.~Ponnusamy and V.~Singh,
Convolution properties of some classes analytic functions,
\textit{Zapiski Nauchnych Seminarov POMI} {\bf 226}(1996), 138--154;
translation in  \textit{J. Math. Sci. (New York)} {\bf 89}(1)(1998), 1008--1020.

\bibitem{RhSch02} Q.~I.~Rahman and G.~Schmeisser,
Analytic Theory of Polynomials, London Mathematical
Society Monographs New Series, 26, Oxford University Press, Oxford, 2002.

\bibitem{rs1} St.~Ruscheweyh and T.~Sheil-Small, Hadamard products of schlicht functions
and the P\'{o}lya-Schoenberg conjecture,
\textit{Comment. Math. Helv.} {\bf 48}(1973), 119--135.

\end{thebibliography}
\end{document}